\documentclass[preprint,12pt]{elsarticle}
\usepackage{graphicx}
\usepackage{subfigure}
\usepackage{epsfig}
\usepackage{amsmath}
\usepackage{amsthm}
\usepackage{amssymb}
\usepackage{lscape}
\usepackage{color}
\usepackage{booktabs}
\usepackage{rotating}
\usepackage{epstopdf,amsfonts,enumitem,color}
\usepackage{float}
\usepackage{hyperref}

\newtheorem{thm}{Theorem}[section]
\theoremstyle{definition}
\newtheorem{example}[thm]{Example}

\numberwithin{equation}{section}
\numberwithin{figure}{section}
\numberwithin{table}{section}

\newcommand{\uu}{{\mathbf u}}

\journal{Computers and Mathematics with Applications}

\begin{document}

\begin{frontmatter}

\title{An experimental comparison of a space-time multigrid method with PFASST for a reaction-diffusion problem}

\author[USI]{Pietro Benedusi\corref{cor1}}
\ead{benedp@usi.ch}
\cortext[cor1]{Corresponding author}

\author[UCB]{Michael L. Minion}
\ead{mlminion@lbl.gov}

\author[USI]{Rolf Krause}
\ead{rolf.krause@usi.ch}

\address[USI]{Institute of Computational Science, University of Italian Switzerland (USI),\\ via Giuseppe Buffi 13, 6900 Lugano, Switzerland}

\address[UCB]{Lawrence Berkeley National Laboratory, Berkeley, CA, USA}


	


\begin{abstract}
We consider two parallel-in-time approaches applied to a (reaction) diffusion problem, possibly non-linear. In particular, we consider PFASST (Parallel Full Approximation Scheme in Space and Time) and space-time multigrid strategies. For both approaches, we start from an integral formulation of the continuous time dependent problem.
Then, a collocation form for PFASST and a discontinuous Galerkin discretization in time for the space-time 
multigrid are employed, resulting in the same discrete solution at the time nodes. Strong and weak scaling of
both multilevel strategies are compared for varying orders of the temporal discretization. Moreover, we investigate 
the respective convergence behavior for non-linear problems  and highlight quantitative differences in execution times.
\end{abstract}



\end{frontmatter}

\section{Introduction} 
Since the clock frequency of computer processors has not increased significantly in the past fifteen years, an increase in computational performance for numerical algorithms can be achieved only by increasing  parallel concurrency, and modern supercomputers now contain many  thousands of computing cores.
Exploiting the capabilities of such massively parallel systems is not straightforward; algorithms with optimal complexity and  excellent scalability must be designed to minimize the run-time of computationally intensive problems, such as
the solution of time dependent partial differential equations (PDEs). When dealing with parallel solvers for discretized PDEs, the solution process is traditionally parallelized in space using domain decomposition techniques, until stagnation. Considering the 
technology trend, the traditional sequential time stepping  will increasingly become the bottleneck for computational scalability for many applications. Hence, the development of new parallel methods that exploit concurrency in the time direction has become essential
for time dependent problems. 
However, parallelization in time can be a challenging task, as, for many physical processes, the time direction is governed by a causality principle, with a preferential direction of information flow through the temporal domain, i.e. forward in time.
Nevertheless, several new methods for temporal parallelization have been proposed in the last 20 years.
For a more comprehensive review regarding the parallel-in-time literature of the past 50 years we refer to \cite{Gander2015_Review}.

The objective of this work is to compare two of the most relevant  recent approaches: PFASST \cite{EmmettMinion2012} and space-time multigrid methods (STMG) \cite{Hackbusch1984,HortonVandewalle1995,FalgoutEtAl2014_MGRIT,NeumuellerGander2016,franco2018multigrid,benedusi2016parallel}.  The current paper is similar in spirit to the comparison presented in \cite{falgout2017multigrid} (where the authors suggest a future comparison with PFASST). Since many parallel-in-time methods are based on a coupling between coarse and fine time propagators, they can be framed in a multilevel-in-time setting; for example, MGRIT \cite{FalgoutEtAl2014_MGRIT}, Parareal in \cite{gander2018multigrid} or PFASST in \cite{bolten2017multigrid}.  Despite the similarities of these different approaches, the methods can behave quite differently on some problems. To date, the majority of  papers on parallel in time methods investigate a single method, and very few have investigated the computational advantages and disadvantages of different methods on well defined benchmarks. Here we attempt to compare to methods using the same spatial discretizations and solvers to emphasis the differences in scaling in the time direction.

In the remainder of this paper we present such a comparison between PFASST and STMG. In Section~\ref{sec_pre} we describe the respective time discretizations of the two approaches. In Section~\ref{sec_problem} we present a reaction-diffusion PDE and its discretization in space and time. In Section~\ref{sec_solvers} we
describe the solution methods
that will be used in Section~\ref{sec_exp}, where weak and strong scaling experiments are reported.

\section{Preliminaries on the time discretizations} \label{sec_pre}
Let us introduce the time discretizations that we use for PFASST and the space-time multigrid (STMG), respectively. Both formulations are equivalent to the same implicit Runge-Kuta (RK) method and are based on an integral form of the continuous problem. To illustrate these methods we consider the initial value problem on a single time step $I_n:=[T_n,T_{n+1}]\subset \mathbb R$
\begin{equation}\label{ode}
u'(t) = f(u(t),t) \quad \text{for} \quad T_n< t < T_{n+1},\qquad u(T_n) = U_0.
\end{equation}

\subsection{Collocation form}
The PFASST algorithm is based on the spectral deferred correction (SDC) method, an iterative scheme introduced in \cite{dutt2000spectral} based on a collocation approximation of \eqref{ode}. Let us consider the Picard integral form of \eqref{ode}
\begin{equation}\label{int_form}
    u(t)=U_0 + \int_{T_n}^t f(u(\tau),\tau)\text{d}\tau,
\end{equation} and the $M$ right Gauss-Radau nodes $\{t_{m}\}_{m=1}^M$ in $I_n$ with $T_n < t_1 < t_2 < \cdots < t_M = T_{n+1} $. We approximate \eqref{int_form} by its collocation form, with $U_m\approx u(t_m)$:
\begin{equation}\label{collocation}
    \mathbf{U} = \mathbf{U}_0 +  \Delta t QF(\mathbf{U}),
\end{equation} where $\Delta t := T_{n+1}-T_n$,
\begin{equation}\label{notation}
\mathbf{U}:=[U_1,...,U_M], \quad \mathbf{U}_0:=[U_0,...,U_0], \quad F(\mathbf{U}):= [f(U_1,t_1),...,f(U_M,t_M)],
\end{equation}
$Q$ is the $M\times M$ matrix $Q:=(q_{m,j})_{m,j=1}^M$ with the quadrature weights
\begin{equation*}
q_{m,j}:=  \frac{1}{\Delta t}\int_{T_n}^{t_m}\ell_j(t)\text{d}t,
\end{equation*}
and $\{\ell_j\}_{j=1}^M$ are the Lagrange polynomials at the $M$ nodes. An SDC iteration can be considered as a preconditioned Richardson iteration to solve \eqref{collocation} (see, e.g. \cite{Huang2006-yo, weiser2015faster}) and, if SDC converges, it is equivalent to an implicit (RK) method, with $q_{m,j}$ being the values in the corresponding Butcher tableaux. The resulting RK method is $A-$stable and has order of accuracy $2M-1$ for $M$ Radau quadrature nodes \cite{Hairer1991-fj}.


\subsection{Discontinuous Galerkin}  \label{DG_section} 
Variational time-stepping methods are receiving increasing interest by the scientific community, especially in the context of adaptivity in space–time, for example in \cite{eriksson1991adaptive,Schmich2008Adaptivity}. Discontinuous Galerkin (DG) methods, in particular, have been widely used to discretize the time direction in the space-time setting as they ensure that the information flows in the positive time direction. They have been employed for a variety of problems such as convection/advection/diffusion equations or the Navier-Stokes equations. For example see the works \cite{jamet1978galerkin,KLAIJ2006589,li1998implementation,SUDIRHAM20061491,Feistauer2011,Besier2012goal,NeumuellerGander2016,benedusi2018space}. The use of DG discretization in time was first introduced in \cite{lasaint1974finite} for the discretization of a neutron transport
equation. In this paper the authors show that, for finite elements of order $q$, the method is strongly $A$-stable, has convergence order $2q+1$ in the nodes, and is equivalent to an implicit (RK) time stepper with $q$ intermediate steps.
The first analysis on DG methods as time stepping techniques is provided by \cite{delfour1981discontinuous} and  \cite{eriksson1985time}, followed by the work of \cite{schieweck2010stable,zhao2014unified,thomee1984galerkin}. More recently, specialized solution methods have been introduced, for example by \cite{smears2016robust,richter2013efficient,hussain2011higher,benedusi2019fast}. A \textit{priori} and \textit{posteriori} error analysis have been also provided, e.g. see \cite{thomee1984galerkin,eriksson1991adaptive,eriksson1995adaptive,Schotzau2010}. See \cite{Shu2014} for a recent survey on the topic. 

Let us consider the weak formulation of \eqref{ode}, where the continuity
at $T_n$ is weakly imposed, and $u\approx U\in \mathbb P_q(I_n)$ (i.e. the space of polynomials of degree $q$) 

\begin{equation}\label{DG_form}
 \int_{T_n}^{T_{n+1}} U'(t) v(t)\,{\rm d}t + (U(T_n) - U_0)v(T_n) = \int_{T_n}^{T_{n+1}} f(U(t),t) v(t)\, {\rm d}t,
\end{equation}for all test functions $v\in \mathbb P_q(I_n)$. Equivalently, integrating by parts \eqref{DG_form}, we obtain the standard DG formulation:

\begin{equation}\label{DG_form2}
 - \int_{T_n}^{T_{n+1}} U(t) v'(t)\,{\rm d}t + U(T_{n+1})v(T_{n+1})- v(T_n)U_0= \int_{T_n}^{T_{n+1}} f(U(t),t) v(t)\,{\rm d}t,
\end{equation} 
where we highlight the \textit{upwind flux} given by the $ v(T_n)U_0$ term. In the interval $I_n$, we construct the approximation $U$ in the nodal form,
\begin{equation} 
  \label{DG_disc}
  U(t) = \sum\limits_{m = 1}^{M} U_m\ell_{n,m}(t),
\end{equation}
where $\left\{\ell_{n,m}\right\}_{m=1}^{M}$ is the basis of Lagrange polynomials of degree $q$ at the $q+1=M$ Gauss-Radau nodes in $I_n$. We can rewrite \eqref{DG_form2}, using the approximation in \eqref{DG_disc} and the definitions in \eqref{notation}, as 
\begin{equation}\label{DG_problem}
    K_q\mathbf{U} = J_q\mathbf{U}_0 + M_qF(\mathbf{U}),
\end{equation} with
\begin{align}
K_q&:=\left[-\int_{T_n}^{T_{n+1}}\ell'_{n,i}(t)\ell_{n,j}(t)\,{\rm d}t + \ell_{n,i}(T_{n+1})\ell_{n,j}(T_{n+1})\right]_{i,j=1}^M, \label{Kq}\\
M_q&:=\left[\int_{T_n}^{T_{n+1}}\ell_{n,i}(t)\ell_{n,j}(t)\,{\rm d}t\right]_{i,j=1}^M,\qquad J_q:=\left[\ell_i(T_n)\ell_j(T_{n+1})\right]_{i,j=1}^M.\label{Mq}
\end{align}

\noindent Let us remark the similarity between \eqref{collocation} and \eqref{DG_problem} and that $J_q\mathbf{U}_0 = [U_0,0,...,0]^T$. For multiple adjacent time elements equation \eqref{DG_problem} can be naturally extended, with obvious notation, as
\begin{equation}\label{DG_problem2}
    K_q\mathbf{U}_n = J_q\mathbf{U}_{n-1} + M_qF(\mathbf{U}_n).
\end{equation}
 

\section{Problem setting and discretization} \label{sec_problem}
Let $\Omega=(0,X)$ be the spatial domain and $T\in\mathbb R^+$ the final time. We consider the following non-linear reaction diffusion equation:
\begin{equation}\label{model}
\left\{\begin{aligned}
&\partial_t u - \partial_{xx} u + \gamma(u^3- u) = 0, && \text{for} \quad(t,x)\in (0,T)\times \Omega,\\
& \partial_x u = 0, && \text{for} \quad(t,x)\in(0, T)\times\partial \Omega, \\
&u=u_0, &&\text{for} \quad t=0 \quad \text{and} \quad x\in\Omega,
\end{aligned}\right.
\end{equation}where $u:=u(t,x)$, $u_0:=u_0(x)$, and $\gamma\geq0$ controls the intensity of the reaction term. Equation \eqref{model} is known as the monodomain model, and it is used to describe the progressive activation of excitable media. For example, in the context of computational medicine, it is employed to simulate the propagation of the electrical potential in the human heart \cite{keener1998mathematical}. The cubic term is a FitzHugh-Nagumo-type reaction, with three zeros $\{-1,0,1\}$ corresponding, respectively, to a resting state, a threshold and an activation state. For $\gamma = 0$ equation \eqref{model} is reduced to the linear heat equation. 

Let $N_t,N_x\in\mathbb N$ be the number of time and space elements respectively, and define the following uniform partitions in time and space:
\begin{equation*}
\begin{aligned}
t_i&:=i\Delta t, &\ \ i&=0,\ldots,N_t, &\ \ \Delta t&:=T/N_t,\\[5pt]
x_j&:=jh, &\ \ j&=0,\ldots,N_x, &\ \ h&:=X/N_x.
\end{aligned}
\end{equation*}
In space, we approximate \eqref{model} with linear finite elements, constructing the discrete operators 
\begin{equation}\label{spatial operators}
K_h:=\left[\int_{\Omega}\varphi_i'(x)\varphi_j'(x){\rm d}x\right]_{i,j=0}^{N_x}, \qquad
M_h:=\left[\int_{\Omega}\varphi_i(x)\varphi_j(x){\rm d}x\right]_{i,j=0}^{N_x},
\end{equation} using the linear Lagrange basis functions  $\{\varphi_i\}_{i=0}^{N_x}\subset H^1(\Omega)$. Referring to Section~\ref{DG_section}, we can consider a space-time finite element approximation of \eqref{model} in $[t_n,t_{n+1}]$ with a tensor structure:
\begin{equation}
    u(x,t) \approx U(x,t)= \sum_{i=1}^{M} \sum_{j=0}^{N_x} u^{n+1}_{i,j}\ell_{n,i}(t)\varphi_j(x)
\end{equation} and assemble the non-linear space-time system of size $N_t M(N_x+1)$ 

\begin{equation}\label{st_system}
\left[\begin{array}{cccc}
A_{q,h} & & & \\
B_{q,h} & A_{q,h} & & \\[4pt]
& \ddots & \ddots & \\[4pt]
& & B_{q,h} & A_{q,h}
\end{array}\right]\left[\begin{array}{c}\uu_1\\ \uu_2\\ \vdots\\ \uu_{N_t}\end{array}\right] + \gamma
(I_{N_t} \otimes M_{q,h})
\left[\begin{array}{c}r(\uu_1)\\ r(\uu_2)\\ \vdots\\ r(\uu_{N_t})\end{array}\right]
=\left[\begin{array}{c}-B_{q,h}\uu_0\\ \mathbf{0}\\ \vdots\\ \mathbf{0}\end{array}\right],
\end{equation} with 
\begin{equation}\label{A_B_ops}
A_{q,h}:= K_q \otimes M_h +  M_q\otimes K_h,\qquad B_{q,h}:= -J_q\otimes M_h,\qquad M_{q,h}:= M_q\otimes M_h,
\end{equation}
$I_{N_t}$ being the identity of size $N_t$. For $n=1,...,N_t$, we have the solution vector
\begin{equation}
    [\uu_n]_{(i-1)N_x+i+j}:=u_{i,j}^n \quad \text{for} \quad i=1,...,M \quad \text{and} \quad j=0,...,N_x,
\end{equation}
with size $M(N_x + 1)$, and the point-wise reaction defined by
\begin{equation}
    [r(\uu_n)]_k:=[\uu_n]_k^3-[\uu_n]_k \quad \text{for} \quad k=0,...,N_xM.
\end{equation}
The initial condition is imposed through $\uu_0:=[0,...,0,u_0(x_0),u_0(x_1),..,u_0(x_{N_x})]$, having $(N_x+1)(M-1)$ zeros. Let use mention that, in the space-time context, increasing $M$ leads to denser blocks. For a detailed description of the weak formulation of \eqref{model} (for $\gamma=0$), the assembly and spectral analysis of system \eqref{st_system}, in a more general finite element framework, we refer to \cite{benedusi2018space}. With respect to standard methods, the storage of system \eqref{st_system} can be expensive in terms of memory; nevertheless, the space-time formulation can be convenient, in terms of scaling and run-time, if \eqref{st_system} is distributed among multiple processors and solved in parallel. Let us remark that we assemble \eqref{st_system} just in the space-time multigrid case; when using PFASST the assembly of the spatial operators in \eqref{spatial operators} is sufficient. For technical limitations related to the current PFASST implementation, we replace the mass matrix $M_h$ with its lumped version for both discretizations.



\section{Solution methods}\label{sec_solvers}
Here we introduce the two solution strategies that will be the object of the comparison.


\subsection{PFASST}\label{sec_pfasst}

\quad The parallel full approximation scheme in space and time (PFASST) was introduced by Emmett and Minion in 2012 \cite{EmmettMinion2012}. As the name suggests,  PFASST can be  described in the context of a multigrid in time method based on a FAS correction on coarse levels \cite{bolten2017multigrid}.
An alternative perspective on how the PFASST method is organized is to view it as a way to perform  SDC iterations for the collocation Eq. \eqref{collocation} on multiple time steps simultaneously.  For parallel efficiency, the SDC iterations are done on a hierarchy of levels as in the multilevel SDC method \cite{SpeckEtAl2015} with communication of new initial conditions passed forward in time between processors after each SDC iteration on each level.  Since the communication is only serial on the coarsest level, the SDC iterations on the finest level are done concurrently, resulting in a potential parallel speedup if the total number of PFASST iterations needed to converge on all the time steps remains relatively small.  One advantage of viewing PFASST from this SDC perspective is  that variants of the original SDC method such as semi-implicit SDC (SISDC) \cite{Minion2004-jn}, can be easily used in the PFASST context.  SISDC methods (also known as implicit-explicit or IMEX) are appropriate for differential equations for which the right hand side of \eqref{ode} can be split into stiff and non-stiff parts.  These methods are often employed in situations where the non stiff component is nonlinear and the stiff term is linear, so that only a linear implicit equation needs to be solved in each time step.  In Section~\ref{sect:nonlinear}, an IMEX treatment is used to treat the nonlinear reaction terms explicitly and the linear diffusion terms implicitly.


\subsection{Space-time multigrid} \label{sec_stmg}
Specialized parallel solvers have been recently developed for large linear systems arising from space-time discretizations. We mention in particular the parallel STMG proposed by \cite{NeumuellerGander2016}, the parallel preconditioners for space-time isogeometric analysis proposed by \cite{hofer2019parallel} and \cite{benedusi2019fast} as well as the block preconditioned GMRES by \cite{mcdonald2016simple}.
When dealing with a space-time discretization, where time is somehow considered as an additional spatial dimension, it is natural to extend the same paradigm for the solving process and consider space-time multigrid type algorithms. 

Multigrid solvers are optimal preconditioners for elliptic problems, and they have proven to be efficient, with some precautions, also for space-time discretizations of parabolic problems. In particular, the heat equation is first order in time and its discretization introduces non symmetric lower bi-diagonal blocks in the space-time system \eqref{st_system}. 

When dealing with anisotropic problems standard multigrid convergence rates deteriorate; see, e.g., \cite{BriHenMcC}. Traditionally there are various ways to address this problem such as accounting for anisotropy in the particular choice of line smoothers and/or adopting a semi-coarsening strategy. In \cite{HortonVandewalle1995,gander2018multigrid,franco2018multigrid,benedusi2020parallel}, for example, the authors explain how the STMG convergence depends critically on the ratio $\mu:=\Delta t/h^2$, unless semi-coarsening strategies are adopted. In particular, for $\mu \ll 1$ (resp. $\mu \gg 1$) coarsening only in time (resp. space) is an effective strategy.

Let us consider a hierarchy of $L$ space-time grids denoted with $l = 1,...,L$ and $l=L$ corresponding to the coarsest one. We construct the space-time restriction operator $I_l^{l+1}$ from level $l$ to level $l+1$ as 
\begin{equation}
I_l^{l+1} = T_l^{l+1} \otimes  M_l^{l+1} \otimes S_l^{l+1}, \quad \text{for} \quad l = 1,...,L-1, \label{space_time_tran}
\end{equation} where $T_l^{l+1}$ and $S_l^{l+1}$ are restriction operators in time and space respectively and $M_l^{l+1}$ is responsible for $M-$coarsening in time, i.e. varying $M$ along the multilevel hierarchy. Definitions of these operators will be provided in the next section. Let us mention that any restriction operator in \eqref{space_time_tran} can be replaced by a suitable identity matrix, resulting in various semi-coarsening strategies. For smoothing, we employ GMRES preconditioned with an incomplete LU factorization (ILU(0)-PGMRES). When multiple parallel cores are used, the space-time system is preconditioned using a block-Jacobi preconditioner, with blocks of size $(N_x+1)N_t M /\text{Cores}$; for each diagonal block the aforementioned PGMRES is then employed. For specialized  preconditioners and corresponding tensor solvers applied to system \eqref{st_system} see, e.g., \cite{pazner2018approximate,borm2001analysis,benedusi2019fast}.
 On the coarsest level, an LU factorization is used and coarse problems are assembled through Galerkin assembly. If $\gamma \neq 0$, equation \eqref{st_system} is non-linear and the STMG algorithm is encapsulated in a Newton iteration. 

\section{Experiments} \label{sec_exp}
\subsection{Implementation}
For the numeric examples in this section, as well as throughout this paper, we used the C++ frameworks PETSc \cite{petsc-user-ref,petsc-web-page} and the embedded domain specific language Utopia\footnote{\url{https://bitbucket.org/zulianp/utopia}} \cite{utopiagit} for the parallel linear algebra and the linear and non-linear solvers. For PFASST we use the modern Fortran library LibPFASST\footnote{\url{https://pfasst.lbl.gov/codes}} that was extended to use the same PETSc data structures and linear solvers as the STMG code to make the comparison as fair as possible. The two discretizations produce, up to machine precision, the same solution in $T$.

Parallel numerical experiments have been performed on the multi-core partition of the  supercomputer Piz Daint of the Swiss national supercomputing centre (CSCS)\footnote{\url{https://www.cscs.ch/computers/piz-daint}}.


\subsection{Solvers specifics and notation}\label{sec_notation}
Next we introduce some of the notation that we are going to use in the following numerical experiments:

\begin{itemize}

\item $\boxed{\text{SMG}^L_{\nu}}$

Multigrid with $L$ levels and spatial coarsening, with $S_l^{l+1}$ in \eqref{space_time_tran} being standard linear bisection, using the stencil $[1\,2\,1]/4$ and the operator $S_l^{l+1}$ having size  $(1+N_x/2^l)\times (1+N_x/2^{l+1})$. As time coarsening is not employed, time transfers in \eqref{space_time_tran} are replaced by identities, i.e. $T_l^{l+1}=I_{N_t}$ and $M_l^{l+1}=I_M$ for all $l=1,...,L-1$. We use V-cycling, with $\nu$ smoothing iterations of ILU(0)-PGMRES.

\item $\boxed{\text{STMG}^L_{\nu}}$

Space-time multigrid with $L$ levels and $T_l^{l+1}$ and $S_l^{l+1}$ in \eqref{space_time_tran} being standard linear bisection operators, both using the same stencil introduced for $\text{SMG}^L_{\nu}$; the operator $T_l^{l+1}$ has size $N_t/2^l\times N_t/2^{l+1}$. We set $M_l^{l+1}=I_M$ for all $l=1,...,L-1$, i.e. $M$ is constant along the multilevel hierarchy. We use V-cycling, with $\nu$ smoothing iterations of ILU(0)-PGMRES.

\item $\boxed{\text{SMMG}^L_{\nu}}$

As $\text{SMG}^L_{\nu}$, but using $M$-coarsening in time; in \eqref{space_time_tran} $M_l^{l+1}$ is obtained through linear interpolation and the number of time nodes $M$ is reduced progressively on the level hierarchy until $M=1$ is reached, i.e. $M=\max\{M-l+1,1\}$ on level $l$.


\vspace{0.1cm}
\item $\boxed{\text{PFASST}^L_{\nu}}$
\vspace{0.1cm}

PFASST solver, as described in Section~\ref{sec_pfasst}, with $L$ levels and $\nu$ sweeps per level. We use ILU(0)-PGMRES as a spatial solution method and standard bisection to create coarse spatial problems.  Regarding temporal coarsening, for performance reasons, we use $M=1$ on all coarse levels.

\end{itemize}
In the numerical results the run-times are expressed in seconds; the assembly of discrete problems and transfer operators are not included in the run-times. The number of iterations to convergence, if present, is reported in square brackets. Convergence is reached when the relative or the absolute preconditioned residual is less then a tolerance of $10^{-9}$. The acronym ``n.c.'' stands for ``not converged'', denoting an increasing residual or if 1000 iterations are exceeded. The tests are restricted to temporal parallelism, i.e. $\#\text{Cores}\leq N_t$, with solvers parameters (i.e. $L$ and $\nu$) which minimize run-time for both approaches. Linear and non-linear iterative solvers are initialized with the zero vector in the space-time case. The spatial diffusion solver in PFASST (PGMRES), are initialized with the best available guess, i.e. the solution at the previous iteration.  

\subsection{Linear example: the heat equation}
In this section we consider the heat equation, i.e. in the following experiments we set $\gamma = 0$ in \eqref{model}, and the initial condition
\begin{equation}\label{init_cond}
u_0(x)=\cos{(\pi x)} + 2\cos{(3\pi x)} + 3\cos{(4\pi x)}\quad \text{for} \quad x\in[0,X],
\end{equation}
with the corresponding analytical solution 
\begin{equation}\label{an_sol}
u(x,t)=\cos{(\pi x)}e^{-\pi^2 t} + 2\cos{(3\pi x)}e^{-9\pi^2 t} + 3\cos{(4\pi x)}e^{-16\pi^2 t}.
\end{equation}
We also consider the analytical solution $\tilde{u}$, obtained after spatial discretization, to focus on the error introduced just by temporal discretization:
\begin{equation}\label{an_sol2}
\tilde{u}(x,t)=\cos{(\pi x)}e^{\rho_1 t} + 2\cos{(3\pi x)}e^{\rho_2 t} + 3\cos{(4\pi x)}e^{\rho_3 t},
\end{equation}
with 
$$ \rho_1 = (2\cos(\pi h)-2)/h^2, \quad \rho_2 = (2\cos(3\pi h)-2)/h^2, \quad \rho_3 = (2\cos(4\pi h)-2)/h^2 .$$

\begin{figure}[H]
    \includegraphics[width=\textwidth]{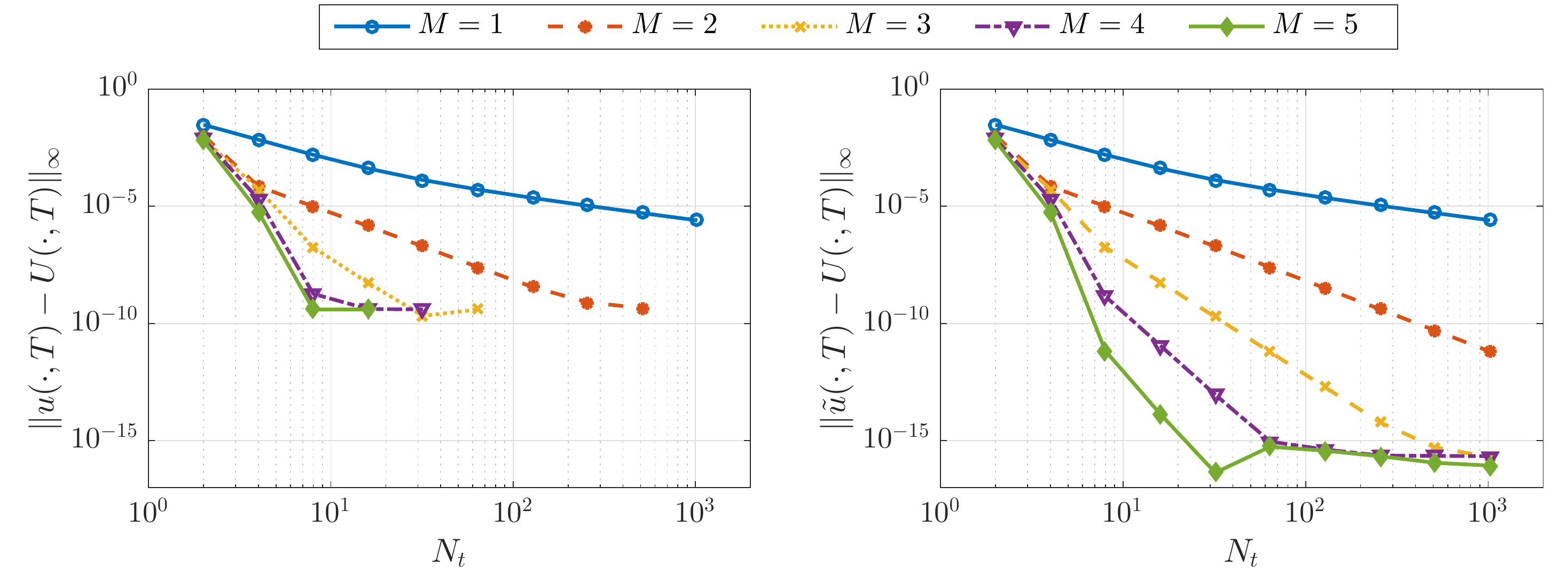}
    \caption{Left: error at the end node w.r.t. the analytical solution \eqref{an_sol}. Right: error w.r.t. \eqref{an_sol2}, i.e. the error of the discrete ODE. 
    In both cases, the errors decrease, with the expected order ($2M-1$), until the spatial error dominates.}
    \label{fig:accuracy}
\end{figure}

\noindent We show, in Figure~\ref{fig:accuracy}, how the error behaves as a  function of the temporal discretization parameters, i.e. $N_t$ and $M$, for problem \eqref{model} with $T=X=1$, $\gamma=0$ and $u_0$ from equation \eqref{init_cond}, discretized, according to \eqref{st_system} with $N_x=1024$. It is possible to observe, in the left plot in Figure~\ref{fig:accuracy}, that the error compared to the exact solution decreases as the number of time steps increases until the spatial error of roughly
$10^{-9}$ dominates. For this reason the solver tolerance is set to $10^{-9}$ in the following numerical experiments. The right-hand plot shows that the temporal error decreases with the correct order $2M-1$ until machine precision is reached.

\noindent In Example~\ref{ex_2}-\ref{ex_3} we report strong and weak scaling results varying the discretization order $M$. We avoid over-resolving in time by reducing $N_t$ as $M$ increases, according to Figure~\ref{fig:accuracy}. In particular, we set $N_t$, depending on $M$, to be the minimum power of two for which $\|u(\cdot,T)-U(\cdot,T)\|_{\infty}<10^{-9}$ is satisfied. With this methodology, the solver tolerance and the spatial and temporal accuracies are approximately the same\footnote{For $M=1$,  we use $N_t=1024$ with a corresponding temporal accuracy of $\sim 10^{-6}$, since $N_t\approx 10^7$ would be required to reach a temporal accuracy of $10^{-9}$.}.

\vspace{0.5cm}

\begin{figure} [H]
    \centering
    \includegraphics[width=\textwidth]{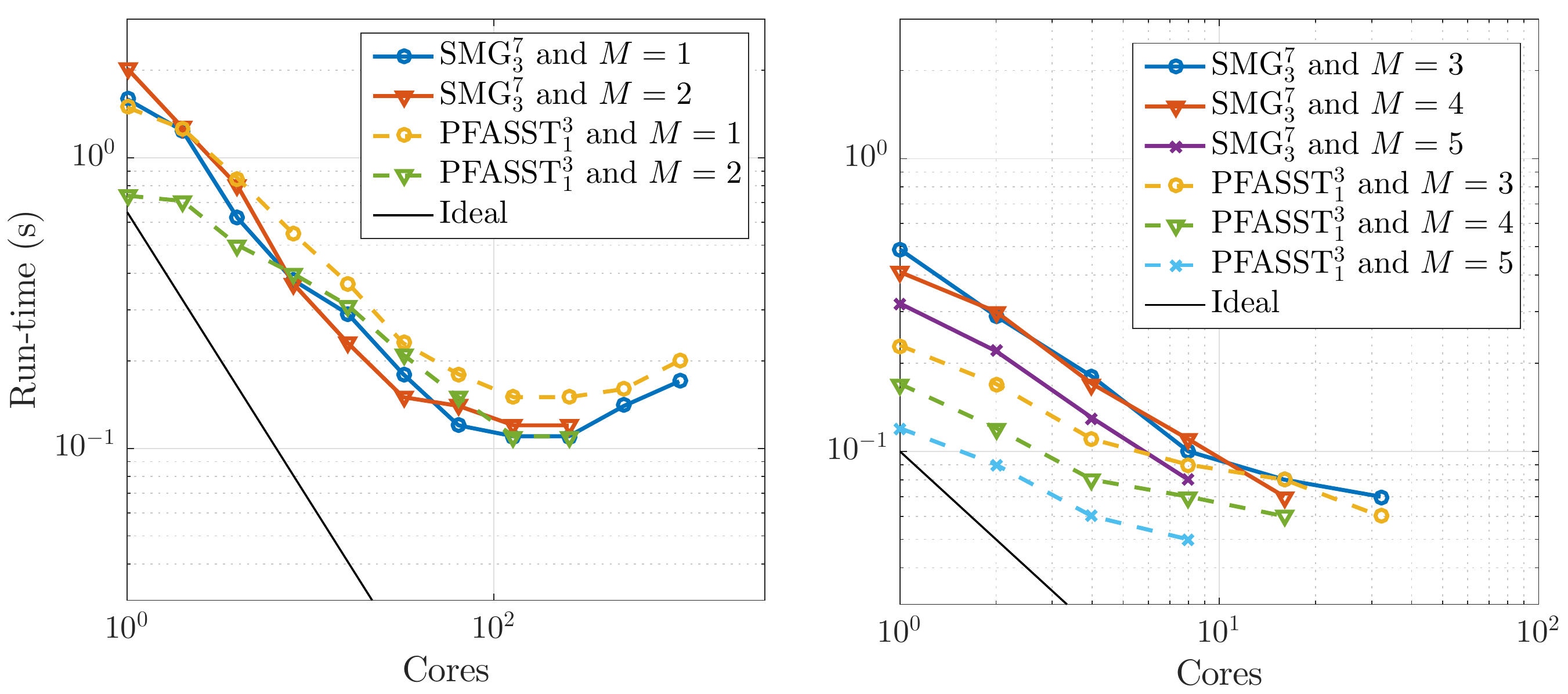}
    \caption{Strong scaling timing results of SMG (solid lines) and PFASST (dashed lines) from Example~\ref{ex_2}; run-times of STMG and SMMG are not included since they are not competitive. We report run-times for $M\in\{1,2\}$ in the left plot and for $M\in\{3,4,5\}$ in the right one.}
    \label{fig:strong}
\end{figure}

\begin{example}\label{ex_2}
(Strong scaling). Let us consider the continuous problem \eqref{model} with parameters $X=T=1$ and $u_0$ from \eqref{init_cond}. We use the discretization parameters $N_x=1024$, $M=\{1,...,5\}$ and $N_t$ varying according to the results of Figure~\ref{fig:accuracy}, to avoid over-resolving in time, except  for $M=1$ where we use $N_t=1024$, with a corresponding accuracy of approximately $10^{-6}$.

We show, in Tables~\ref{tab_ex21}--\ref{tab_ex23}, run-times and iterations of the three multilevel approaches described in Section~\ref{sec_notation}. The most relevant results are illustrated in Figure~\ref{fig:strong}.

\begin{table}[H]
\centering
{\small 
\begin{tabular}{l|l|l|l|l|l}
&  \multicolumn{5}{c}{SMG$_3^7$} \\[2pt]
\hline
     & $M=1,$ & $M=2,$ & $M=3,$ & $M=4,$ & $M=5,$  \\[2pt]
    Cores & $N_t=1024$ &$N_t=256$ & $N_t=32$ &$N_t=16$ &$N_t=8$  \\[2pt]
    \hline\hline
    1  & 1.59 [2] & 2.03 [4] & 0.49 [4] & 0.41 [5] & 0.32 [5] \\[2pt]
    2  & 1.24 [2] & 1.27 [4] & 0.29 [4] & 0.30 [5] & 0.22 [5]  \\[2pt]
    4  & 0.62 [2] & 0.80 [4] & 0.18 [4] & 0.17 [5] & 0.13 [5]  \\[2pt]
    8  & 0.38 [2] & 0.37 [4] & 0.10 [4] & 0.11 [5] & 0.08 [5]  \\[2pt]
    16 & 0.29 [3] & 0.23 [4] & 0.08 [4] & 0.07 [5] \\[2pt]
    32 & 0.18 [3] & 0.15 [4] & 0.07 [4] \\[2pt]
    64 & 0.12 [3] & 0.14 [4]  \\[2pt]
   128 & 0.11 [3] & 0.12 [4]  \\[2pt]
  256  & 0.11 [3] & 0.12 [4]  \\[2pt]
  512  & 0.14 [3] \\[2pt]
  1024 & 0.17 [3] \\[2pt]
 \end{tabular}} 
 \caption{Seven level space-time multigrid run-times and iterations to convergence with no temporal coarsening.} \label{tab_ex22}
 \end{table}
 
 \begin{table}[H]
\centering
{\small 
\begin{tabular}{l|l|l|l|l|l}
&  \multicolumn{5}{c}{STMG$_3^5$} \\[2pt]
\hline

     & $M=1,$ & $M=2,$ & $M=3,$ & $M=4,$ & $M=5,$  \\[2pt]
    Cores & $N_t=1024$ &$N_t=256$ & $N_t=32$ &$N_t=16$ &$N_t=8$  \\[2pt]
    \hline\hline
    1  & 2.66 [6] & 6.57 [18] & 2.12 [30] & 1.75 [30] & 1.23 [30] \\[2pt]
    2  & 1.53 [6] & 3.71 [18] & 1.21 [29] & 0.98 [30] & 0.71 [31] \\[2pt]
    4  & 0.80 [6] & 1.82 [18] & 0.64 [32] & 0.64 [33] & 0.43 [32] \\[2pt]
    8  & 0.50 [6] & 1.05 [18] & 0.39 [34] & 0.52 [37] & 0.27 [34] \\[2pt]
    16 & 0.35 [6] & 0.67 [18] & 0.29 [37] & 0.26 [37] \\[2pt]
    32 & 0.22 [6] & 0.53 [18] & 0.22 [38] \\[2pt]
    64 & 0.17 [6] & 0.45 [17] \\[2pt]
   128 & 0.17 [6] & 0.37 [17] \\[2pt]
  256  & 0.17 [6] & 0.37 [17] \\[2pt]
  512  & 0.18 [7] \\[2pt]
  1024 & 0.29 [9] \\[2pt]
  
  
 \end{tabular}} 
 \caption{Five level space-time multigrid run-times and iterations to convergence, with full space-time coarsening.} \label{tab_ex21}
 
 \end{table}

 \begin{table}[H]
\centering
{\small 
\begin{tabular}{l|l|l|l|l}
&  \multicolumn{4}{c}{SMMG$_3^7$} \\[2pt]
\hline
     & $M=2,$ & $M=3,$ & $M=4,$ & $M=5,$  \\[2pt]
    Cores & $N_t=256$ & $N_t=32$ &$N_t=16$ &$N_t=8$  \\[2pt]
    \hline\hline
    1  & 3.01 [12] & 0.99 [17] & 0.79 [15] & 0.75 [18] \\[2pt]
    2  & 1.77 [12] & 0.62 [17] & 0.44 [15] & 0.50 [18] \\[2pt]
    4  & 0.86 [12] & 0.28 [17] & 0.23 [15] & 0.28 [18]  \\[2pt]
    8  & 0.47 [12] & 0.17 [17] & 0.15 [15] & 0.15 [18]  \\[2pt]
    16 & 0.24 [12] & 0.12 [17] & 0.11 [15] \\[2pt]
    32 & 0.15 [12] & 0.09 [17] \\[2pt]
    64 & 0.12 [12]   \\[2pt]
    128 & 0.11 [12]  \\[2pt]
   256  & 0.11 [12]    \\[2pt]
  
 \end{tabular}} 
 \caption{Seven level space-time multigrid run-times and iterations to convergence, with $M-$coarsening in time. The column for $M=1$ is not present as it would be equivalent to the one of Table~\ref{tab_ex22}.} \label{tab_ex23}
 
 \end{table}

\begin{table}[H]
\centering
{\small 
\begin{tabular}{l|c|c|c|c|c}
&  \multicolumn{5}{c}{PFASST$_1^3$} \\[2pt]
\hline
  & $M=1,$ & $M=2,$ & $M=3,$ & $M=4,$ & $M=5,$  \\[2pt]
    Cores & $N_t=1024$ &$N_t=256$ & $N_t=32$ &$N_t=16$ &$N_t=8$  \\[2pt]
    \hline
    1  & 1.49 [2]  & 0.74 [3] &  0.23 [11] & 0.17 [12] & 0.12 [12] \\[2pt]
    2  & 1.26 [3]  & 0.71 [11] & 0.17 [14] &  0.12 [14] &  0.09 [15]  \\[2pt]
    4  & 0.85 [4]  & 0.50 [13] & 0.11 [16] & 0.08 [17] & 0.06 [18]  \\[2pt]
    8  & 0.55 [6]  & 0.40 [17] & 0.09 [20] & 0.07 [21] & 0.05 [22]  \\[2pt]
    16 & 0.37 [7]  & 0.31 [24] & 0.08 [28] &  0.06 [29] \\[2pt]
    32 & 0.23 [7]  & 0.21 [29] & 0.06 [35] \\[2pt]
    64 & 0.18 [7] & 0.15 [30]  \\[2pt]
    128 & 0.15 [8] & 0.11 [30]  \\[2pt]
  256  & 0.15 [7] & 0.11 [30]  \\[2pt]
  512  & 0.16 [8] \\[2pt]
  1024 & 0.20 [7] \\[2pt]
 \end{tabular}} 
 \caption{Three level PFASST run-times and iterations to convergence.} \label{tab_ex24}
 \end{table}


\end{example}

\newpage
The results of Example~\ref{ex_2} suggest several observations.
First, for the space-time multigrid methods, using spatial coarsening only (SMG)  obtains the fastest space-time multigrid convergence and run-times. We can explain this behavior using the discretization parameter
\begin{equation*}
    \mu=\frac{\Delta t}{\Delta x^2}=\frac{N_x^2}{N_t}.
\end{equation*}
In all cases considered $\mu\gg 1$; from space-time multigrid literature (e.g. \cite{franco2018multigrid},\cite{horton1995space}) we can expect time coarsening to be not effective in this scenario, as we observe by the larger iteration counts in Table~\ref{tab_ex21}. The case of $\mu\leq1$ would not be meaningful in this setting, as it would result in a unnecessary over-resolved time discretization, with $N_t\geq N_x^2>10^6$. On the other hand, the results collected in Table~\ref{tab_ex23} suggest that $M-$coarsening is the most convenient way to coarsen in time in the space-time multigrid case, since run-times are close to the case of spatial coarsening only. Note also that the 
$M-$coarsening in time for SMMG and PFASST are the same. Comparing SMG with PFASST, cf. Figure~\ref{fig:strong}, we observe an overall similar scaling and run-times, especially for $M\in\{1,2\}$. As expected from the space-time paradigm, SMG is characterized by a slower sequential run-time compensated by a better scaling. As as results, both algorithm  achieve a similar stagnation run-time.  For higher values of $M$, PFASST is somewhat faster  than SMG, presumably  
due to the reduced cost of the coarsest level.  Note also that all the methods achieve some parallelism even using 2 processors. Finally, it is important to note that for the same accuracy, the higher-order methods  have lower run times than the lower-order methods, often with fewer processors.  Hence, although the scaling in terms of iterations is better for the lower-order methods, they are more expensive in practice.

\begin{figure}[H]
    \centering
    \includegraphics[width=\textwidth]{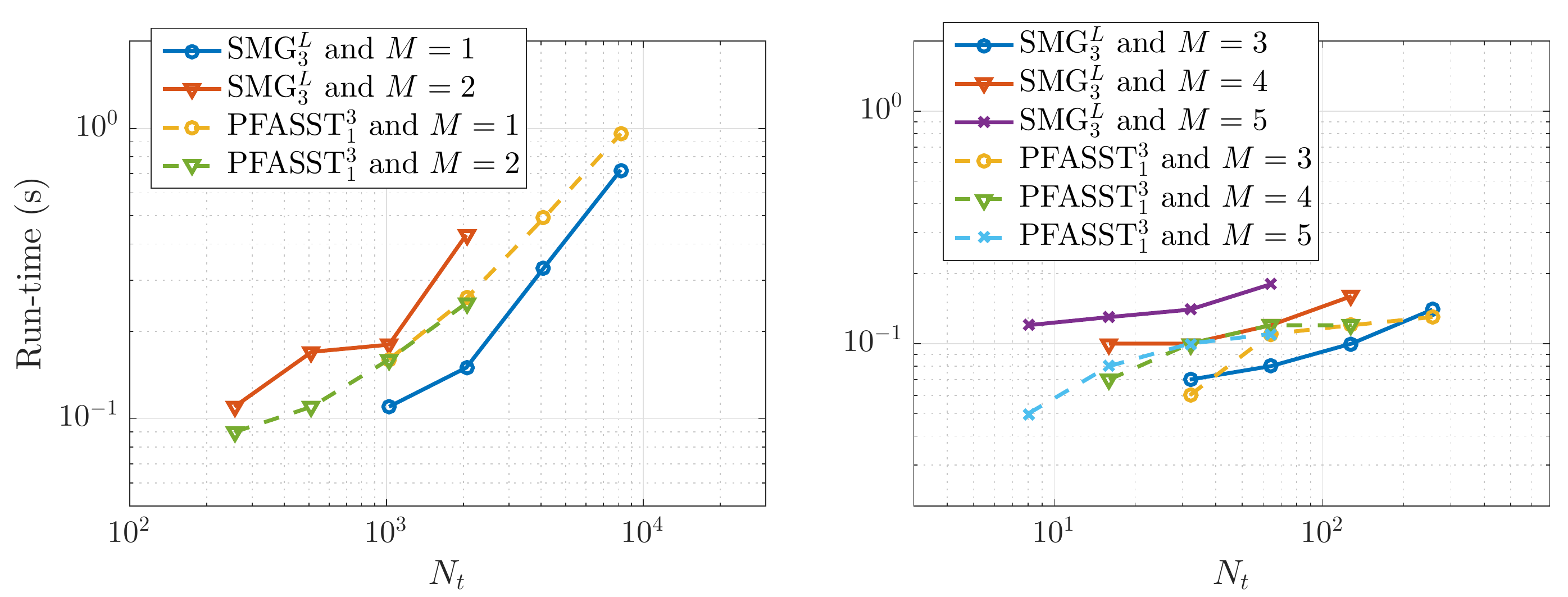}
    \caption{Weak scaling timing results of SMG (solid lines) and PFASST (dashed lines) from Tables~\ref{tab_ex31}--\ref{tab_ex32}; run-times of STMG and SMMG are not included since they are not competitive. We report run-times for $M\in\{1,2\}$ in the left plot and for $M\in\{3,4,5\}$ in the right one.}
    \label{fig:weak}
\end{figure}

\begin{example} \label{ex_3}
(Weak time scaling in $N_t$) Let us consider the continuous problem \eqref{model} with parameters $X=T=1$ and $u_0$ from \eqref{init_cond}. We use the discretization parameters $N_x=1024$, $M=\{1,...,5\}$ and $N_t=C_M\cdot \text{Cores}$. The parameter $C_M$ depends on $M$ and is chosen according to the parallel saturation from Example~\ref{ex_2}; for example, for $M=1$ and $N_t=1024$, according to Tables~\ref{tab_ex22}--\ref{tab_ex24}, we have maximum speedup with 256 cores and therefore $C_1=4$. Similarly $C_2=2$ and $C_M=1$ for $M\geq3$. We point out that the accuracy of the solution as a function of $N_t$ is varying in the tables according to Figure~\ref{fig:accuracy}: doubling $N_t$ corresponds to a higher accuracy as $M$ increases. 
We report, in Tables~\ref{tab_ex31}--\ref{tab_ex32} and Figure~\ref{fig:weak}  run-times and iterations of the multilevel approaches described in Section~\ref{sec_notation}. 
 \end{example}
 
 \begin{table}[H]
    \begin{minipage}{.5\linewidth}
      \centering
        {\small 
\begin{tabular}{r|r|c|c|c}
\multicolumn{5}{c}{$M=1$} \\[2pt]
    \hline  
     Cores & $N_t$ & L & time [its.] & $R$\\[2pt]
    \hline \hline
      256  & 1024 & 7 & 0.11 [3] & 1  \\[2pt]
      512  & 2048 & 8 & 0.15 [4] & 1.4  \\[2pt]
      1024 & 4096 & 8 & 0.33 [5] & 3.0\\[2pt]
      2048 & 8192 & 8 & 0.72 [6] & 6.5
  
 \end{tabular}
 } 
    \end{minipage}%
    \begin{minipage}{.5\linewidth}
    \qquad
      {\small 
\begin{tabular}{r|r|c|c|c}
\multicolumn{5}{c}{$M=2$} \\[2pt]
    \hline  
     Cores & $N_t$ & L & time [its.] & $R$\\[2pt]
    \hline \hline
     128  & 256 & 7 & 0.11 [4] & 1 \\[2pt]
     256  & 512 & 7 & 0.17 [4] & 1.5\\[2pt]
     512  & 1024 & 8 & 0.18 [4] & 1.6\\[2pt]
     1024 & 2048 & 8 & 0.43 [5] & 3.9 
  
 \end{tabular}
 } 
    \end{minipage} 
    
    \vspace{6ex}
    
    \begin{minipage}{.5\linewidth}
      \centering
        {\small 
\begin{tabular}{r|r|c|c|c}
\multicolumn{5}{c}{$M=3$} \\[2pt]
    \hline  
    Cores & $N_t$ & L & time [its.] & $R$\\[2pt]
    \hline \hline
     32  & 32 & 7 & 0.07 [4] & 1 \\[2pt]
     64  & 64 & 7 & 0.08 [5] & 1.1 \\[2pt]
     128  & 128 & 8 & 0.10 [5] & 1.4\\[2pt]
     256  & 256 & 8 & 0.14 [5] & 2.0 
  
 \end{tabular}
 } 
    \end{minipage}%
    \begin{minipage}{.5\linewidth}
    \qquad
      {\small 
\begin{tabular}{r|r|c|c|c}
\multicolumn{5}{c}{$M=4$} \\[2pt]
    \hline  
    Cores & $N_t$ & L & time [its.] & $R$\\[2pt]
    \hline \hline
    16  & 16 & 7 & 0.10 [5] & 1\\[2pt]
    32  & 32 & 7 & 0.10 [5] & 1.0\\[2pt]
    64  & 64 & 7 & 0.12 [5]  & 1.2\\[2pt]
    128 & 128 & 8 & 0.16 [6]  & 1.6 
  
 \end{tabular}
 } 
    \end{minipage} 
    
    \vspace{6ex}
    \centering
          {\small 
\begin{tabular}{r|r|c|c|c}
\multicolumn{5}{c}{$M=5$} \\[2pt]
    \hline  
    Cores & $N_t$ & L & time [its.] & $R$\\[2pt]
    \hline \hline
     8   & 8 & 7 & 0.12 [5] & 1 \\[2pt]
     16  & 16 & 7 & 0.13 [6] & 1.1  \\[2pt]
     32  & 32 & 7 & 0.14 [7] & 1.2\\[2pt]
     64  & 64 & 7 & 0.18 [7]  & 1.5 
  
 \end{tabular}
 } 
    \vspace{1.0cm}
    \caption{Weak scaling in time of a seven level space-time multigrid SMG$_L^3$, with no temporal coarsening. The ratio $R$ is computed dividing the current run-time by the base one (in the first line for each table) and $R=1$ denotes an ideal weak scaling. Since no temporal coarsening is present, the weak scaling is poor for $M\in\{1,2\}$.} \label{tab_ex31}
\end{table}



 \begin{table}[H]
    \begin{minipage}{.5\linewidth}
      \centering
        {\small 
\begin{tabular}{c|r|r|c}
\multicolumn{4}{c}{$M=1$} \\[2pt]
    \hline  
    Cores & $N_t$ & time [its.] & $R$\\[2pt]
    \hline \hline
     256  & 1024 & 0.16 [7] & 1.0 \\[2pt]
     512  & 2048 & 0.26 [8] & 1.6 \\[2pt]
     1024  & 4096 & 0.49 [10] & 3.1 \\[2pt]
     2048  & 8192 & 0.96 [10] & 6.0
  
 \end{tabular}
 } 
    \end{minipage}%
    \begin{minipage}{.5\linewidth}
    \qquad
      {\small 
\begin{tabular}{c|r|r|c}
\multicolumn{4}{c}{$M=2$} \\[2pt]
    \hline  
    Cores & $N_t$ & time [its.] & $R$\\[2pt]
    \hline \hline
     128  & 256 & 0.09 [30] & 1.0\\[2pt]
     256  & 512 & 0.11 [30] & 1.2 \\[2pt]
     512  & 1024 & 0.16 [29] & 1.8 \\[2pt]
     1024  & 2048 & 0.25 [30] & 2.8
  
 \end{tabular}
 } 
    \end{minipage} 
    
    \vspace{3ex}
    
    \begin{minipage}{.5\linewidth}
      \centering
        {\small 
\begin{tabular}{c|r|r|c}
\multicolumn{4}{c}{$M=3$} \\[2pt]
    \hline  
    Cores & $N_t$ & time [its.] & $R$\\[2pt]
    \hline \hline
     32  & 32 & 0.06 [35] & 1.0 \\[2pt]
     64  & 64 & 0.11 [38] & 1.8  \\[2pt]
     128  & 128 & 0.12 [38] & 2.0 \\[2pt]
     256  & 256 & 0.13 [36] & 2.2
  
 \end{tabular}
 } 
    \end{minipage}%
    \begin{minipage}{.5\linewidth}
    \qquad
      {\small 
\begin{tabular}{c|r|r|c}
\multicolumn{4}{c}{$M=4$} \\[2pt]
    \hline  
    Cores & $N_t$ & time [its.] & $R$\\[2pt]
    \hline \hline
     16  & 16 & 0.07 [29] & 1.0 \\[2pt]
     32  & 32 & 0.10 [38] & 1.4\\[2pt]
     64  & 64 & 0.12 [42] & 1.4\\[2pt]
     128  & 128 & 0.12 [42] & 1.7
  
 \end{tabular}
 } 
    \end{minipage} 
    
    \vspace{3ex}
    \centering
          {\small 
\begin{tabular}{c|r|r|c}
\multicolumn{4}{c}{$M=5$} \\[2pt]
    \hline  
    Cores & $N_t$ & time [its.] & $R$\\[2pt]
    \hline \hline
     8  & 8 &  0.05 [22] & 1.0\\[2pt]
     16  & 16 & 0.08 [30] & 1.6 \\[2pt]
     32  & 32 & 0.10 [39] & 2.0\\[2pt]
     64  & 64 & 0.11 [42] & 2.2
  
 \end{tabular}
 } 
    \vspace{0.7cm}
    \caption{Weak scaling in time of PFASST$^3_1$. The ratio $R$ is computed dividing the current run-time by the base one (in the first line for each table) and $R=1$ denotes an ideal weak scaling. We can notice that the weak scaling for $M=1$ is poor, since no temporal coarsening is present in this case.} \label{tab_ex32}
\end{table}

It is clear from Figure~\ref{fig:weak}
that higher-order methods display better weak scaling than the lower-order ones. This is partly due to the limited time-coarsening in both cases, that would become relevant as $N_t$ grows. It should be noted that for $M>1$, the error w.r.t. the analytical solution is already saturated in the first rows of Tables~\ref{tab_ex31}--\ref{tab_ex32} (according to Figure~\ref{fig:accuracy}) and increasing $N_t$ does not produce a more accurate solution. In this scaling regime, there is no consistent winner between SMG and PFASST for all $M$ and $N_t$.

\subsection{Non-linear example: the monodomain equation.}\label{sect:nonlinear}
In this section we consider the full reaction-diffusion model, i.e. $\gamma >0$ in \eqref{model}. Note that for the space-time discretization in \eqref{st_system}, a non-linear solver is required due to the cubic reaction term. For the implementation in PFASST an IMEX or semi-implicit  method is employed \cite{Minion2003-qr} treating the  
nonlinear reaction terms explicitly.  Hence the cost per iteration of the PFASST method is essentially the same as for the linear case since the implicit part is much more expensive.

To model a traveling wave in an excitable media we consider reaction dominated examples. In this case we set a narrow initial stimulus in the centre of the domain and we chose $T$ such that the final solution is stationary, i.e. for all $x$ we have $u(T,x) \simeq 1$ and $\partial_t u(T,x)\simeq 0$. 

\vspace{0.5cm}

\begin{example} \label{ex_4}
(Strong scaling)  Let us consider the continuous problem \eqref{model} with model parameters $X=10,T=2,\gamma = 5$ and the initial condition 
\begin{equation*}
u_0 = 2\exp{\left(\frac{x-X}{0.1}\right)^2}.    
\end{equation*}
We use the discretization parameters $N_x=N_t=1024$ and $M=\{1,...,5\}$. We show, in Table~\ref{tab_ex41}, run-times and Newton iterations of the space-time strategy, using SMG$_3^7$ as linear solver and, in Table~\ref{tab_ex42}, the PFASST data. 
\end{example}

Note that in the cases where PFASST converges, the run time is significantly smaller than the corresponding SMG times.  As mentioned above, this is due to the fact that the PFASST implementation is using a semi-implicit or IMEX time stepping method, while  SMG method is fully implicit requiring Newton iterations.  The failure of PFASST to converge for large time steps is also due to the IMEX stepping, which has a time step restriction due to the explicit treatment of the reaction term (see, e.g. \cite{Minion2004-jn}). 
 The reaction term could also be handled implicitly in PFASST using a multi-implicit approach \cite{Bourlioux2003-bq} as was done in \cite{Gotschel2019-ka} to increase the stability, but we defer this sort of comparison to future work.

\begin{table}[H]
\centering
{\small 
\begin{tabular}{l|l|l|l|l|l}
&  \multicolumn{5}{c}{SMG$^7_3$} \\[2pt]
\hline
     & $M=1,$ & $M=2,$ & $M=3,$ & $M=4,$ & $M=5,$  \\[2pt]
    Cores & $N_t=1024$ &$N_t=256$ & $N_t=32$ &$N_t=16$ &$N_t=8$  \\[2pt]
    \hline\hline
    1  & 46.4 [16]  & 31.1 [16] & 6.47 [16] & 5.18 [16] & 301 [58] \\[2pt]
    2  & 29.7 [16]  & 19.5 [16] & 4.50 [16] & 3.75 [16] & 123 [45]  \\[2pt]
    4  & 16.6 [16]  &  11.0 [16] & 2.61 [16] & 2.17 [16] & 125 [47]  \\[2pt]
    8  & 10.6 [16]  & 7.45 [16] &  1.81 [16] & 1.60 [16] & 34.4 [18]  \\[2pt]
    16 & 8.51 [16]  & 5.81 [16] & 1.36 [16]. & 1.16 [16] \\[2pt]
    32 & 7.22 [16]  & 4.92 [16] & 1.07 [16] \\[2pt]
    64 & 6.11 [16] & 4.54 [16]  \\[2pt]
    128 & 5.10 [16] & 3.90 [16]  \\[2pt]
  256 & 5.50 [16] & 5.48 [16]  \\[2pt]
  512 & 5.64 [16] \\[2pt]
  1024 & 6.20 [16] \\[2pt] 
 \end{tabular}} 
 \caption{Run-time of a seven level space-time multigrid, with no temporal coarsening and corresponding Newton iterations.} \label{tab_ex41}
 \end{table}
 
\begin{table}[H]
\centering
{\small 
\begin{tabular}{l|l|l|l|l|l}
&  \multicolumn{5}{c}{PFASST$_1^3$} \\[2pt] 
\hline
     & $M=1,$ & $M=2,$ & $M=3,$ & $M=4,$ & $M=5,$  \\[2pt]
    Cores & $N_t=1024$ &$N_t=256$ & $N_t=32$ &$N_t=16$ &$N_t=8$  \\[2pt]
    \hline\hline
    1  & 1.51 [2]  & 0.54 [5] & 0.23 [93] & n.c. & n.c. \\[2pt]
    2  & 1.16 [3]  & 0.42 [6] & 0.20 [94] & n.c. & n.c.  \\[2pt]
    4  & 0.67 [4]  & 0.25 [8] & n.c. & n.c. & n.c. \\[2pt]
    8  & 0.40 [5]  & 0.16 [8] & n.c. & n.c. & n.c.  \\[2pt]
    16 & 0.27 [6]  & 0.11 [8] & n.c. & n.c. \\[2pt]
    32 & 0.18 [6]  & 0.07 [8] & n.c. \\[2pt]
    64 & 0.14 [6] & 0.07 [8]  \\[2pt]
    128 & 0.13 [6] & 0.07 [8]  \\[2pt]
  256  & 0.15 [6] & 0.08 [8]  \\[2pt]
  512  & 0.16 [6] \\[2pt]
  1024 & 0.20 [6] \\[2pt]
 \end{tabular}} 
 \caption{Three level PFASST run-times and iterations, n.c. abbreviating “not converged”.} \label{tab_ex42}
 \end{table} 

\section{Conclusions}  
In this paper  we discussed the parallel performance of multilevel space-time solution strategies and of the algorithm PFASST for a (reaction) diffusion problem. 
From an implementation prospective, space-time multigrid approaches are convenient since the time parallelization boils down to the parallel solution of a system of equations (in \eqref{st_system}), for example using fast and parallel preconditioned Krylov methods as PGMRES. A tensor structure between space and time grids allows for a flexible choice of coarsening strategies, since transfer operators in \eqref{space_time_tran} can be set independently. On the other hand, the assembly of system \eqref{st_system} comes at a cost, in terms of time and, especially, memory footprint. Such cost can be be reduced significantly when \eqref{st_system} is distributed among many processors and if highly parallel assembly routines are used, as parallel Kronecker products in \eqref{A_B_ops}. 


In Examples~\ref{ex_2}--\ref{ex_3} we investigated the scalability of different parallel iterative strategies for a diffusion problem. We obtained similar performance from PFASST and the parallel space-time multigrid with no temporal coarsening (SMG). The use of high order methods in time, reducing the number of time steps $N_t$ accordingly, is convenient for both approaches, in terms of overall performance and especially for the weak scaling in time.

As expected from the literature, full space-time coarsening or time coarsening are not effective in the settings we considered ($\mu\gg 1$). In the space-time multigrid framework $M-$coarsening in time can be advantageous w.r.t. coarsening in the number of time steps $N_t$, in terms of stability, but employing just coarsening in space remains the best option for the  discretizations considered.

In Example~\ref{ex_4} we considered a non-linear reaction-diffusion problem. For such a problem the space-time approach is limited to a fully implicit treatment of the non-linearity and the corresponding use of a non-linear solver, such as  Newton's method. In particular, we observe that the number of Newton iterations to convergence is not robust in terms of problem parameters and initial guess. On the other hand, in this respect PFASST is more flexible since it allows one to treat the non-linearity explicitly, through an IMEX approach. Such a strategy, even if less stable (especially for large $\Delta t$ and high order $M$), can reduce dramatically the time-to-solution.

\section*{Acknowledgements}
\quad The authors acknowledge the Deutsche Forschungsgemeinschaft (DFG) as part of the ``ExaSolvers'' Project in the Priority Programme 1648 ``Software for Exascale Computing'' (SPPEXA) and the Swiss National Science Foundation (SNSF) under the lead agency grant agreement SNSF-162199. The work of M. Minion was supported by the U.S. Department of Energy, Office of Science,
Office of Advanced Scientific Computing Research, Applied Mathematics program
under contract number DE-AC02005CH11231. This research used resources of the National Energy Research
Scientific Computing Center, a DOE Office of Science User Facility
supported by the Office of Science of the U.S. Department of Energy
under Contract No. DE-AC0205CH11231.

 \textbf{Declaration of interest:} none.

\bibliographystyle{elsarticle-num}
\bibliography{main}

\end{document}